\documentclass[11pt]{amsart}

\usepackage{fullpage,amsbsy,amssymb,enumitem}

\title[Modelling Silicosis: Structure of Equilibria]{Modelling
  Silicosis: The Structure of Equilibria}

\author{Fernando P. da Costa}

\address[F.P. da Costa]{Univ. Aberta, Dep. of Sciences and Technology,
  Rua da Escola Polit\'ecnica 141-7, P-1269-001 Lisboa, Portugal, and
  Univ. Lisboa, Instituto Superior T\'ecnico, Centre for Mathematical
  Analysis, Geometry and Dynamical Systems, Av. Rovisco Pais,
  P-1049-001 Lisboa, Portugal.}  \email{fcosta@uab.pt}

\thanks{Research partially supported by FCT/Portugal through project
  UID/MAT/04459/2013.}

\author{Michael Drmota} \address[M. Drmota]{TU Wien, Institute of
  Discrete Mathematics and Geometry, Wiedner Hauptstrasse 8-10, A-1040
  Vienna, Austria.}  \email{michael.drmota@tuwien.ac.at.}

\thanks{Research partially supported by the Austrian Science
  Foundation FWF, project F 50-02.}

\author{Michael Grinfeld} \address[M. Grinfeld]{Univ. of Strathclyde,
  Dep. of Mathematics and Statistics, 26 Richmond Street, Glasgow G1
  1XH, United Kingdom.}  \email{m.grinfeld@strath.ac.uk.}

\date{January 24, 2019}

\subjclass{Primary 34A34; Secondary 34E05, 92C45}

\keywords{Coagulation--fragmentation--death equations, silicosis,
  asymptotics, Mellin transform}

\newcommand{\F}{\mathcal{F}}
\newcommand{\K}{\mathcal{K}}
\newcommand{\be}{\begin{equation}}
\newcommand{\ee}{\end{equation}}
\newcommand{\Nb}{{\mathbb N}}
\newcommand{\Rb}{{\mathbb R}}

\renewcommand{\H}{\mathcal{H}}

\newcommand{\blot}[1]{}

\DeclareMathOperator{\sgn}{sgn}

\newtheorem{theo}{Theorem} 
\newtheorem{lemma}[theo]{Lemma}
\newtheorem{cor}[theo]{Corollary}
\newtheorem{prop}[theo]{Proposition}

\newenvironment{Proof}{\removelastskip \vskip12pt plus 1pt \noindent
	{\em Proof.\/}\rm }{\hfill$\square$ \vskip12pt plus 1pt}

\begin{document}

\begin{abstract}
  We analyse the structure of equilibria of a
  coagulation--fragmentation--death model of silicosis. We present exact
  multiplicity results in the particular case of piecewise-constant
  coefficients, results on existence and non-existence of
  equilibria in the general case, as well as precise asymptotics for
  the infinite series that arise in the case of power law
  coefficients.  
\end{abstract}

\maketitle

\section{Introduction}

We examine the equations considered in \cite{Tran} for the dynamics of
alveolar macrophages faced with an inhalation of quartz particles in
the lungs. The equations are of coagulation type, though not
understood in the standard way. Thus the interest of the problem is
triple: it is of medical and environmental health importance to
understand how the system reacts to continuous exposure to quartz; the
model shows how versatile coagulation equations are; and the
mathematical structure of the system (which we do not a priori
truncate as is done in \cite{Tran}) presents interesting challenges.
In this paper we only discuss the model itself and the structure of
its equilibria, leaving proofs of global existence and stabilisation
results to future work.

If we denote by $M_i$ the concentration of macrophages containing $i$
quartz particles (which we will call the $i$-th cohort), by $x$ the
concentration of quartz, by $r$ (which can be a function of $x$ the
rate of supply of new macrophages, following \cite{Tran} we obtain the
following equations:
\begin{equation}\label{eq1}
  \begin{array}{l}
    \dfrac{dM_0}{dt}  ~=~ r - k_0 x M_0 - (p_0+q_0) M_0,\\  
		\mbox{} \\
    \dfrac{dM_i}{dt}  ~=~ k_{i-1} x M_{i-1} - k_i x M_i -(p_i+q_i)
    M_i, \;\; i \geq 1,
  \end{array} 
\end{equation}
where $k_i$ is the rate of phagocytosis of a macrophage containing
$i$ particles of quartz, $p_i$, is the transfer rate of macrophages in
the $i$-th cohort to the muco-ciliary escalator, i.e. the rate of
their removal together with their quartz baggage, and $q_i$ is the
rate of death of the $i$-th cohort which results in the release of the
quartz burden.  

Note that unlike \cite{Tran} we do not impose an upper limit on  the
number of particles a macrophage can contain. What is not done in
\cite{Tran} is to provide an equation for the evolution of the
concentration of $x$; their interest is in the system dynamics
following an instance of inhalation, while we are more concerned with
analysing system behaviour under continuous influx of quartz. Thus we
add to (\ref{eq1}) the following equation:
\begin{equation}\label{eq2}
  \frac{dx}{dt} = \alpha - x \sum_{i=0}^\infty  k_i M_i +
  \sum_{i=0}^\infty q_i i M_i.
\end{equation}  
Here $\alpha$ the rate of inhalation of quartz.

Thus the object of our study is the system (\ref{eq1})--(\ref{eq2}),
considered as an infinite-dimensional dynamical system on a suitable
sequence space. Before we analyse it, let us remark that hence
(\ref{eq1})--(\ref{eq2}) is an example of a coagulation-death system,
in which the ``monomers'' (quartz particles) are structurally
different from ``clusters'' (cells containing these particles); this
shows the versatility of coagulation-fragmentation framework, and in
particular its suitability to describe phagocytosis phenomena (e.g. of
neutrophils consuming bacteria).

As in \cite{Tran} we make the assumptions that $k_i$ and $p_i$ are
non-increasing in $i$. We allow $q_i$ to grow with $i$.

The model of \cite{Tran} is biologically sophisticated, also involving
neutrophils and communication between neutrophils and macrophages.
In (\ref{eq1}), $r$  should express the amount of ``distress'' in the
system, embodied in the number of macrophages with more than a
sublethal load of quartz, i.e. those that are more likely to die and
release their load than to be removed via the muco-ciliatory
escalator.  In other words, if the sublethal load is $s$ particles per
cell, a biologically reasonable assumption is that
$r$ is a bounded increasing function of $\sum_{i=s+1}^\infty M_i$ (see
eqs. 7--8 in \cite{Tran}). In the present work we take $r$ to be a
constant, but our analysis here illuminates the more general case
described above as well.

A simple instance of allowable coefficients for which the structure of
equilibria can be analysed explicitly will be considered in section~\ref{constantcoeff} below. 
The structure of the equilibria in a more general case, where the coefficients satisfy some power 
law relations, will be considered in section~\ref{powercoeff}.  

\section{Equilibria}\label{equil}

We start by proceeding formally and then justify our steps in the sections below. 
Suppose system \eqref{eq1}-\eqref{eq2} has an equilibrium. Then the $M_0$
equation at equilibrium can be solved for $M_0$ in terms of $x$ (and
$r$) to give
\[
  M_0 = \frac{r}{k_0 x + p_0+q_0}.
\]
Similarly, $M_1$ will be given by
\[
  M_1 = \frac{rk_0x}{(k_0 x+ p_0+q_0) (k_1 x + p_1+q_1)}.
\]
Continuing recursively, we have
\[
  M_i = \frac{r x^i \prod_{j=0}^{i-1} k_j}{\prod_{j=0}^i (k_j x+ p_j +
    q_j)}.
\]
Setting $d_i= (p_j+q_j)/k_j$, we can rewrite this as
\begin{equation}\label{mi}
 M_i = \frac{r x^i}{k_i\prod_{j=0}^i (x+ d_j)},\qquad i \geq 0.
\end{equation} 

\subsection{A piecewise-constant class of coefficients}\label{constantcoeff}

A simple instance of allowable coefficients is to take all $k_i$ equal to $k$
and,
\[
p_i = \begin{cases} 1 & \text{if $i \leq N,$}\\ 0 & \text{if $i \geq N+1,$}\end{cases}
\quad\text{and}\quad
q_i = \begin{cases} 0 & \text{if $i \leq N,$}\\ 1 & \text{if $i \geq N+1.$} \end{cases}
\]

Then $d_j = 1/k$, and using \eqref{mi} we easily compute 
\be
\sum_{i=0}^\infty k_iM_i = k\sum_{i=0}^\infty M_i = \frac{r}{x+1/k}\sum_{i=0}^\infty\left(\frac{x}{x+\/k}\right)^i = rk, \label{serie1} 
\ee
and
\begin{eqnarray}
\sum_{i=0}^\infty iq_iM_i & = & \sum_{i = N+1}^\infty iM_i \nonumber \\
& = & \frac{1}{k}\frac{r}{x+1/k}\sum_{i = N+1}^\infty i\left(\frac{x}{x+1/k}\right)^i \nonumber \\
& = & \frac{1}{k}\frac{r}{x+1/k}\left(\frac{x}{x+1/k}\right)^{N+1}
      \sum_{i = 0}^\infty (i+N+1)\left(\frac{x}{x+1/k}\right)^i \nonumber \\
& = & \frac{1}{k}\frac{r}{x+1/k}\left(\frac{x}{x+1/k}\right)^{N+1}\left(\sum_{i = 0}^\infty i\left(\frac{x}{x+1/k}\right)^i +
      (N+1)\sum_{i = 0}^\infty \left(\frac{x}{x+1/k}\right)^i\right) \nonumber \\
& = & \frac{1}{k}\frac{r}{x+1/k}\left(\frac{x}{x+1/k}\right)^{N+1}\Bigl(k^2x(x+1/k) + k(N+1)(x+1/k)\Bigr) \nonumber \\
& = & r\left(\frac{x}{x+1/k}\right)^{N+1}\Bigl(kx + (N+1)\Bigr).\label{serie2}
\end{eqnarray}
Thus, plugging \eqref{serie1} and \eqref{serie2} into equation for the equilibrium quartz concentration we obtain
\be
\frac{\alpha}{r} - \mathcal{F}_{N,k}(x) = 0 \label{eqF}
\ee
where
\be
\mathcal{F}_{N,k}(x) := kx\left(1-\left(\frac{x}{x+1/k}\right)^{N+1}\right) - (N+1)\left(\frac{x}{x+1/k}\right)^{N+1}. \label{defF}
\ee

\medskip

\begin{prop}\label{th2}
	
	For all $r, k >0$ and $N\in \Nb,$ there is $\alpha^*$ such that \eqref{eqF} has no solutions if $\alpha>\alpha^*.$
\end{prop}

\begin{Proof}
It suffices to observe that $\mathcal{F}_{N,k}(0)=0$, $\mathcal{F}'_{N,k}(0)=k>0$ and
\begin{eqnarray*}
\lim_{x\to +\infty}\mathcal{F}_{N,k}(x) 
& = & \lim_{x\to +\infty}\left[kx\frac{(x+1/k)^{N+1}-x^{N+1}}{(x+1/k)^{N+1}}
-(N+1)\left(\frac{x}{x+1/k}\right)^{N+1}\right]  \\
& = & \lim_{x\to +\infty}\left[(N+1)\left(\frac{x}{x+1/k}\right)^{N+1} + O(x^{-1}) - (N+1)\left(\frac{x}{x+1/k}\right)^{N+1}\right] \\
& = & 0. 
\end{eqnarray*} 
This implies that $\mathcal{F}_{N,k}$ has an absolute maximum in $\Rb^+.$ 
Defining $\alpha^* := \displaystyle{r\max_{\Rb^+}\mathcal{F}_{N,k}},$
the result follows.
\end{Proof}

\medskip

We now prove that, for each $\alpha \in (0, \alpha^*)$, there are
exactly two solutions of \eqref{eqF}.

\medskip

\begin{prop}\label{th3}
	
	Let $r, k >0$ and $N\in \Nb.$ Let $\alpha^* := \displaystyle{r\max_{\Rb^+}\mathcal{F}_{N,k}}.$  Then, for every 
	$\alpha \in (0, \alpha^*)$ there are exactly two solutions of \eqref{eqF}.
	
\end{prop}

\begin{Proof}
	To prove the result we establish that 
	$\mathcal{F}_{N,k}$ has a single stationary point in
	$\Rb^+$, which, then, must be the absolute maximum whose existence was 
	established above. This, together with the already proved facts that
	 $\mathcal{F}_{N,k}(0) = 0$ and 
	$\lim_{x\to +\infty}\mathcal{F}_{N,k}(x) =0$, proves the result.
	
	Let $y:= \frac{x}{x+1/k}$. Then $x= \frac{y/k}{1-y}$ and
	\begin{eqnarray}
	\widetilde{\mathcal{F}}_{N,k}(y) := \mathcal{F}_{N,k}(x(y)) 
	& = & \frac{y}{1-y}\left(1-y^{N+1}\right) -(N+1)y^{N+1} \nonumber \\
	& = & \frac{y}{1-y}\left(1-y^{N+1}-(N+1)(1-y)y^{N}\right) \nonumber \\
          & = & \frac{y}{1-y}\underbrace{\left(1-(N+1)y^N+Ny^{N+1}
                \right)}_{=:f_N(y)}\nonumber
	\end{eqnarray}
	Since $\frac{dx}{dy} = \frac{1}{k}\frac{1}{(1-y)^2}>0$, we have 
	$\sgn \widetilde{\mathcal{F}}'_{N,k}(y) = \sgn \mathcal{F}'_{N,k}(x(y))$. Thus we need only to study the
	function in the new variable $y\in [0,1)$. 
	Observing that  $f'_N(y)=-N(N+1)y^{N-1}+N(N+1)y^N = -N(N+1)y^{N-1}(1-y),$ we
	have 
	\begin{eqnarray}
	\widetilde{\mathcal{F}}'_{N,k}(y) 
	& = & \frac{1}{(1-y)^2}f_N(y) + \frac{y}{1-y}f_N'(y) \nonumber \\
	& = & \frac{1}{(1-y)^2}\left(1-(N+1)y^N +Ny^{N+1}-N(N+1)(1-y)^2y^N\right)\nonumber \\
	& = & \frac{1}{(1-y)^2}\underbrace{\left(1-(N+1)^2y^N+N(2N+3)y^{N+1}-N(N+1)y^{N+2}\right)}_{=:p_N(y)}. \nonumber
	\end{eqnarray}
	Let us consider the polynomial $p_N$ in $[0,1].$ It is clear
        that $p_N(0)=1$ and $p_N(1)=0.$ Its derivative is
        $p_N'(y) = N(N+1)y^{N-1}q_N(y)$, where
        $q_N(y):= -(N+1) + (2N+3)y - (N+2)y^2.$ We easily conclude
        that the zeros of $q_N(y)$ are $y_1 = \frac{N+1}{N+2}$ and
        $y_2=1,$ and that
        $\sgn\left(y-\frac{N+1}{N+2}\right) \sgn q_N(y) > 0.$ This
        means that $p_N$ has a minimum at $y=\frac{N+1}{N+2}$ and must
        be an increasing function in the interval
        $\left(\frac{N+1}{N+2}, 1\right)$. Since $p_N(1)=0,$ this
        implies the value of $p_N(y)$ at $y=\frac{N+1}{N+2}$ must be
        negative, which, together with $p_N(0)>0$ and the fact that
        $p_N$ is strictly decreasing in
        $\left(0,\frac{N+1}{N+2}\right)$, means, due to the
        intermediate value theorem, that there is one, and only one,
        zero of $p_N$ in this set, and hence in $(0,1)$, i.e., there
        is a single stationary point of $\mathcal{F}_{N,k}$ in
        $\Rb^+$.
\end{Proof}

\medskip

\subsection{Power type coefficients}\label{powercoeff}

We consider now the more complex case of coefficients satisfying
some power relations.

\begin{theo}\label{th1}
Let $M_i$ be given by (\ref{mi}). Assume that $z=\inf_i d_i >0$. 
Assume also that $q_i/k_i$
grows no faster than a power of $i$. Then for all $x \geq 0$
\[
  \sum_{i=0}^\infty k_i M_i < \infty \qquad \text{and} \qquad
  \sum_{i=0}^\infty i q_i M_i < \infty.
\]
\end{theo}

\begin{Proof}
This follows by the Ratio Test, as
\[
\frac{k_{i+1} M_{i+1}}{k_i M_i} = \frac{x}{x+d_{i+1}} \leq \frac{x}{x+z}<1,
\]
for all $i$. Also,  
\[
\frac{(i+1)q_{i+1} M_{i+1}}{i q_i M_i} = \frac{1+i}{i}
\frac{(q_{i+1}/k_{i+1})}{q_i/k_i} \frac{x}{x+d_{i+1}}.
\]
Pick $\epsilon = {\displaystyle \frac12 \left( \frac{z+x}{x} -1
  \right)}$. We can find $N=N(x)$ such that
\[
\frac{1+i}{i}
\frac{(q_{i+1}/k_{i+1})}{q_i/k_i} \leq (1+ \epsilon)
\]
for all $i \geq N$. But then for all $i\geq N$ we have that
\[   
\frac{(i+1)q_{i+1} M_{i+1}}{i q_i M_i} \leq \frac12 \left( \frac{x}{x+z} +1
  \right) < 1.
\]
\end{Proof}

So the equation for equilibrium quartz concentration can be written in the form
\begin{equation}\label{eqequil}
\frac{\alpha}{r} = \frac{x}{x+d_0} + \frac{1}{x+d_0}
\sum_{i=1}^\infty \left( x - i \frac{q_i}{k_i} \right)
\prod_{j=1}^i \frac{x}{x+d_j}
 =:
\F(x). 
\end{equation}

Our first main result provides a quite general sufficient condition for
the existence of equilibria.

\begin{theo}\label{thmain1}
Assume that $z:=\inf_i d_i >0$ and 
also  that $q_i/k_i$
grows no faster than a power of $i$. Let $\rho_i := p_i/k_i.$ 

If $d_i = o(i \rho_i)$ (as $i\to\infty$) then  $\F(x) \to \infty$ (as $x\to\infty$)
and consequently we have an equilibrium (\ref{eqequil}) for all $\alpha,r$.

If $i\rho_i = O(d_i)$ (as $i\to\infty$) then $\F(x)$ is bounded. Thus
there exists $m> 0$ such that there exists an equilibrium for
$\alpha/r < m$ and no equilibrium for $\alpha/r  \geq m$.

Finally if $i\rho_i = o(d_i)$ (as $i\to\infty$) then $\F(x)$ is
bounded and we have $\F(x) \to 0$ (as $x\to\infty$). In this case
there exists $m> 0$ such that there exists an equilibrium for
$\alpha/r \le m$ and no equilibrium for $\alpha/r > m$.
\end{theo}

In order to prove Theorem~\ref{thmain1} we have to study $\F(x)$ in more detail.

Since $d_i := \frac{p_i+q_i}{k_i} = \rho_i + \frac{q_i}{k_i}$, we have
\begin{eqnarray*}
\F(x) &~=~&   \frac{x}{x+d_0} + \frac{1}{x+d_0}\sum_{i=1}^\infty 
\Bigl(x- i\frac{q_i}{k_i}\Bigr)\prod_{j=1}^i\frac{x}{x+d_j} \nonumber \\
&~=~&  \frac{x}{x+d_0} + \frac{1}{x+d_0}\sum_{i=1}^\infty 
\Bigl(x- id_i + i\rho_i\Bigr)\prod_{j=1}^i\frac{x}{x+d_j} \nonumber \\
&~=~&  \underbrace{\left\{\frac{x}{x+d_0} - \frac{1}{x+d_0}\sum_{i=1}^\infty 
\Bigl(id_i-x\Bigr)\prod_{j=1}^i\frac{x}{x+d_j}\right\}}_{=: \mathcal{G}(x)} +
\frac{1}{x+d_0}\sum_{i=1}^\infty 
i\rho_i\prod_{j=1}^i\frac{x}{x+d_j} \label{Fx2}
\end{eqnarray*}

\begin{prop}\label{th4}
With the above assumptions and notation, we have that $\mathcal{G}(x)=
0,$ $\forall x\geq 0$.
\end{prop}

\begin{Proof}
Since $\mathcal{G}(0) = 0$ it is sufficient to consider the case $x> 0$.
We prove that the sum $S$ of the series
	\begin{equation}
	\frac {1}{x} \sum_{i=1}^\infty \Bigl(id_i-x\Bigr)\prod_{j=1}^i\frac{x}{x+d_j} \label{serie}
	\end{equation}
	is equal to $1$ for all values of $x>0.$ Let $S_n$ denote the partial sums of \eqref{serie} and set 
\[
a_n = 1 - S_n.
\]
We will show by induction that
\begin{equation}\label{eqan}
a_n =  \frac{(n+1) x^n }{\prod_{j=1}^n (x+d_j) }.
\end{equation}
Obviously we have for $n=1$
\[
a_1 ~=~ 1 - S_1 = 1- \frac{1}{x}\sum_{i=1}^1
\left(id_i-x\right)\prod_{j=1}^i\frac{x}{x+d_j}
1- \frac{1}{x}(d_1-x)\frac{x}{x+d_1} = \frac{2x}{x+d_1},
\]
as required. Assume (\ref{eqan}) is true for some $n$.
Then 
\begin{eqnarray*}
a_{n+1} &~=~& 1 - S_{n+1} \\ 
&~=~& 1- \Bigl(S_n +
\frac{1}{x}\left((n+1)d_{n+1}-x\right)\prod_{j=1}^{n+1}
\frac{x}{x+d_j}\Bigr) \\ 
&~=~& a_n -
x^{n}\left((n+1)d_{n+1}-x\right)
\prod_{j=1}^{n+1}\frac{1}{x+d_j}\\ 
&~=~&
\frac{(n+1) x^n }{\prod_{j=1}^n (x+d_j) } - 
x^{n}\left((n+1)d_{n+1}-x\right)
\prod_{j=1}^{n+1}\frac{1}{x+d_j}\\ 
&~=~&
x^n \frac{ (n+1)(x+d_{n+1}) - ((n+1)d_{n+1}-x)}
{ \prod_{j=1}^{n+1}(x+d_j)  }\\ 
&~=~&
\frac{(n+2) x^{n+1} }{\prod_{j=1}^{n+1} (x+d_j) }
.
\end{eqnarray*}	
This proves (\ref{eqan}) for all $n\ge 1$.

Now $a_n \to 0$ as $n\to\infty$ follows trivially from
\[
0 < a_n = (n+1)\frac{x^n}{\displaystyle{\prod_{j=1}^n(x+d_j)}} <
(n+1)\Bigl(\frac{x}{x+z}\Bigr)^n \longrightarrow 0,\hbox{ as } n\to
  +\infty,
\]
So we conclude that $S=\lim S_n = \lim (1-a_n) = 1$. Hence
$\mathcal{G}(x)=0$ holds also for all $x>0.$
\end{Proof}

Using Proposition~\ref{th4} we conclude that, with the power law
assumptions on the coefficients, $\F$ can be written as
\begin{equation}
\F(x) = \frac{1}{x+d_0}\sum_{i=1}^\infty
i \rho_i \prod_{j=1}^i\frac{x}{x+d_j} =: \frac{1}{x+d_0}\H(x), \label{newF}
\end{equation}

With the help of the next proposition we can get some information 
on the growth order of $\H(x)$.

\begin{prop}\label{th4.2}
Assume that $z=\inf_i d_i >0$. Then we have, for $x\ge 0$,
\[
\sum_{i=1}^\infty d_i \prod_{j=1}^i  \frac x{x+d_j} = x.
\]
\end{prop}

\begin{Proof}
The equality is trivially satisfied for $x=0$. Thus we just have
to consider the case $x> 0$, where we set
\[
b_n = x - \sum_{i=1}^n d_i \prod_{j=1}^i  \frac x{x+d_j}.
\]
We prove by induction that 
\begin{equation}\label{eqbn}
b_n = \frac{x^{n+1}}{\prod_{j=1}^n (x + d_j)}.
\end{equation}
This is clearly true for $n=1$:
\[
b_1 = x - d_1 \frac x{x+d_1} = \frac {x^2} {x+d_1}.
\]
Now assume that (\ref{eqbn}) is satisfied for some $n\ge 1$.
Then we have
\begin{eqnarray*}
b_{n+1} &=& b_n - d_{n+1} \prod_{j=1}^{n+1}  \frac x{x+d_j} \\
&=& \frac{x^{n+1}}{\prod_{j=1}^n (x + d_j)} - \frac{d_{n+1} x^{n+1}} {\prod_{j=1}^{n+1} (x + d_j)} \\
&=& x^{n+1} \frac{ x + d_{n+1} - d_{n+1}  } {\prod_{j=1}^{n+1} (x + d_j)} \\
&=& \frac{x^{n+2}}{\prod_{j=1}^{n+1} (x + d_j)}
\end{eqnarray*}
as proposed. 
Finally we have
\[
0 < b_n \le x \left( \frac x{x+z} \right)^n \to 0, \quad \mbox{as $n\to\infty$}.
\]
This implies $\lim b_n = 0$ and proves the proposition for $x> 0$.
\end{Proof}	

We are now ready to prove Theorem~\ref{thmain1}.

\begin{Proof}
First suppose that $d_i = o(i\rho_i)$ (as $i\to\infty$). 
Fix some $\varepsilon> 0$ and suppose that $d_i \le \varepsilon i \rho_i$
for $i\ge i_0 = i_0(\varepsilon)$. Then we have (also by applying 
Proposition~\ref{th4.2})
\begin{eqnarray*}
\H(x) &=& \sum_{i=1}^\infty i\rho_i \prod_{j=1}^i \frac x{x+d_j} \\
&=& \sum_{i=1}^{i_0-1} i\rho_i \prod_{j=1}^i \frac x{x+d_j} + \sum_{i=i_0}^\infty i\rho_i \prod_{j=1}^i \frac x{x+d_j} \\
&\ge& \sum_{i=1}^{i_0-1} i\rho_i \prod_{j=1}^i \frac x{x+d_j} + \frac 1{\varepsilon} \sum_{i=i_0}^\infty d_i \prod_{j=1}^i \frac x{x+d_j} \\
&=& \sum_{i=1}^{i_0-1} \left( i\rho_i - \frac 1\varepsilon d_i\right) \prod_{j=1}^i \frac x{x+d_j} 
+ \frac 1{\varepsilon} \sum_{i=1}^\infty d_i \prod_{j=1}^i \frac x{x+d_j} \\
&=& O(1) + \frac x{\varepsilon}.
\end{eqnarray*}
Consequently 
\[
\liminf_{x\to\infty} \frac {\H(x)}x \ge \frac 1\varepsilon.
\]
Since $\varepsilon > 0$ can be arbitrarily chosen we have $\H(x)/x \to \infty$ (as $x\to\infty$) and, thus,
\[
\F(x) = \frac {\H(x)}{x+d_0} \to \infty, \quad \mbox {as $x\to\infty$}.
\]
Since $\F(0) = 0$ and $\F(x)$ is continuous it follows that there exists an equilibrium in all cases.

Next suppose that $i\rho_i = O(d_i)$, that is, there is a constant $K> 0$ such that $i\rho_i \le K d_i$ for all $i\ge 1$.
Hence,
\[
\H(x) = \sum_{i=1}^\infty i\rho_i \prod_{j=1}^i \frac x{x+d_j} \le K\sum_{i=1}^\infty d_i \prod_{j=1}^i \frac x{x+d_j} = K\, x
\]
and consequently $\F(x)$ is bounded. Clearly, if we set
\[
m:= \sup_{x\ge 0} \F(x)
\]
then there exists no equilibrium if $\alpha/r > m$ and again since $\F(0) = 0$ and by continuity 
there is an equilibrium if $\alpha/r < m$.

Finally if $i\rho_i = o(d_i)$ (as $i\to\infty$) then is follows (as above) that $F(x) = o(x)$ (as $x\to\infty$).
By continuity there exists
\[
m:= \max_{x\ge 0} \F(x).
\]
Hence, then there exists no equilibrium if $\alpha/r > m$ and an equilibrium if $\alpha/r \le m$.
\end{Proof}

Consider now the case where the coefficients are given by the
following power laws:
\[
p_i = i^{-p},
\qquad
q_i = i^q,
\quad
\hbox{ and }
\quad
k_i = i^{-k},
\]    
for $i\in\Nb^+$ and nonnegative constants $p, q$ and $k$, Let $p_0,
q_0$ and $k_0\neq 0$ be given. Then, writing $a:=q+k\geq 0$ and
$b:=k-p \in \Rb,$ we have $d_i =i^{a} +i^{b}$ and $\rho_i = i^{b}.$

By a direct application of Theorem~\ref{thmain1} we get the following property:

\begin{cor}\label{Th2}
Suppose that $a\ge 0$ and $b > -2$.

If $b>a-1$, (\ref{eq1})--(\ref{eq2}) has 
an equilibrium for all $\alpha/r$. 

If $b = a-1$, there is a value $m>0$ such that for $\alpha/r<m$,
(\ref{eq1})--(\ref{eq2}) has an equilibrium 
and no equilibria if $\alpha/r\geq m$.

If $b< a-1$, there is a value $m>0$ such that for $\alpha/r \leq m$,
(\ref{eq1})--(\ref{eq2}) has  an equilibrium
and if $\alpha/r>m$, there are none.
\end{cor}

\section{Precise Asymptotics}

It is also an interesting problem to obtain precise asymptotics
for the case where $d_i =i^{a} +i^{b}$ and $\rho_i = i^{b}$, $b<a$. 
In order to make our analysis slightly easier we will concentrate on the
case
\[
d_i =i^{a} \quad \mbox{and}\quad  \rho_i = i^{b}.
\]
We will use the following notation:
\begin{equation}\label{eqK}
  \K_{a,b}(x) = \sum_{i=1}^\infty
  i^{b+1}\prod_{j=1}^i\frac{x}{x+j^a}, \; \; \;
  \H_{a,b}(x) = \sum_{i=1}^\infty
  i^{b+1}\prod_{j=1}^i\frac{x}{x+j^a+j^b}.
\end{equation}

Then we have the following result:

\begin{theo}\label{Th1}
Suppose that $a>0$, $b > -2$ and let $\K_{a,b}(x)$ be given by (\ref{eqK}). 
Then as $x\to\infty$, $\K_{a,b}(x)$ admits an asymptotic expansion
such that 
\[
\K_{a,b}(x) \sim \frac{\Gamma\left( \frac{b+2}{a+1} \right)}{
(a+1)^{1-(b+2)/(a+1)}} x^{\frac{b+2}{a+1}} + O \left( x^{\frac{b+1}{a+1}}
\right) 
\]
if $b > -1$ and  
   
\[
\K_{a,b}(x) \sim \frac{\Gamma\left( \frac{b+2}{a+1} \right)}{
(a+1)^{1-(b+2)/(a+1)}} x^{\frac{b+2}{a+1}} + O \left( 1 \right)  
\]
if $b \leq -1$.
\end{theo}

The proof is mainly based on the following asymptotic series
representation:
\begin{lemma}\label{LemmaR}
Let $R(a,A;v)$ denote the infinite sums
\[
R(a,A;v) = \sum_{i=1}^\infty  e^{-v i^{a+1}} i^A,
\]
where $a, v> 0$, and $A$ is real. The following holds:
\begin{itemize}
\item[1.]
If $(A+1)/(a+1)$ is different from $0,-1,-2,\ldots$, then as $v\to 0$, 
\begin{equation}\label{Ranonint}
  R(a,A;v) \sim \frac {\Gamma\left( \frac {A+1}{a+1} \right) }{a+1}
  v^{-\frac {A+1}{a+1} } + \sum_{k=0}^\infty \frac{(-1)^k}{k!}
  \zeta(-A-k(a+1)) v^{k},
\end{equation}
where $\zeta(\cdot)$ is the Riemann zeta function.
\item[2.]
If $(A+1)/(a+1) = -k_0$ for some integer $k_0\ge 0$ then, as $v\to 0$,
\begin{equation}\label{Raint}
    R(a,A;v)  \sim  \frac{(-1)^{k_0}}{(a+1)k_0!}\left( 1 + (H_{k_0} +
      a\gamma)\log \frac 1v \right) v^{k_0} 
    + \sum_{k\ge 0,\, k\ne k_0} \frac{(-1)^k}{k!} \zeta(-A-k(a+1))
    v^{k},
\end{equation}
where $H_k = 1 + \frac 12 + \cdots \frac 1k$ denotes the $k$-th
harmonic number, $H_0=0$, and $\gamma = 0.5772156\ldots$ is the
Euler--Mascheroni constant.
\end{itemize}
\end{lemma}

\begin{Proof}
  We recall (see, e.g., \cite[Part I]{FGD}) that the Mellin transform
  of a function $f(v)$ is given by
\[
\hat f(s) = \int_0^\infty f(v) v^{s-1}\, dv,
\]
and converges usually in a strip $a_1 < \Re(s) < a_2$. Under suitable
regularity assumptions (for example that $f(v)$ is continuous and of
bounded variation) the function $f(v)$ can be recovered from the
integral
\[
  f(v) = \frac 1{2\pi i} \lim_{T\to\infty} \int_{C-iT}^{C+iT} \hat
  f(s) v^{-s}\ ds,
\]
where $a_1 < C < a_2$. 

In our case it is an easy exercise to show that the Mellin transform
of $R(a,A;v)$ is given by
\begin{equation}\label{eqintrep}
  \hat R(a,A;s) = \int_0^\infty R(a,A;v) v^{s-1}\, dv = \Gamma(s)
  \zeta((a+1)s - A),
\end{equation}
The integral converges for $\Re(s) > \max\bigl(0,\frac{A+1}{a+1}\bigr)$.
Consequently we have
\[
  R(a,A;v) = \frac 1{2\pi i} \int_{C-i\infty}^{C+i\infty} \Gamma(s)
  \zeta((a+1)s - A) v^{-s}\, ds,
\]
where $C> \max\bigl(0,\frac{A+1}{a+1}\bigr)$. Since the $\Gamma$-function
decreases exponentially fast on vertical lines we could replace the
limit $\lim_{T\to\infty} \int_{C-iT}^{C+iT}$ by the indefinite
integral $\int_{C-i\infty}^{C+i\infty}$.

The function $\Gamma(s) \zeta((a+1)s - A)$ has a meromorphic
continuation to the whole complex plane.  The only singularities are
simple poles coming from $\Gamma(s)$ at $s = -k$ with residue
\[
{\rm Res}(\Gamma(s),-k) = \frac{(-1)^k}{k!} \qquad (k = 0,1,2,\ldots).
\]
and the simple pole of $\zeta((a+1)s - A)$ at $s = \frac{A+1}{a+1}$
with residue
\[
{\rm Res}\left(\zeta((a+1)s - A),\frac{A+1}{a+1} \right) = \frac 1{a+1}.
\]
The idea is to shift the integral in (\ref{eqintrep}) to the left and
to collect residues of the polar singularities that are passed. There
are (again) no convergence problems of the integral due to the
$\Gamma$ factor.

Assume first that $\frac{A+1}{a+1}$ is different from
$0,-1,-2,\ldots$. If we shift the integral to
$\Re(s) = -M- \frac 12$, where $M$ is a positive integer and
$-M < \frac{A+1}{a+1}$ then we have
\begin{eqnarray*}
  \lefteqn{R(a,A;v) =}\\ &=& \frac {\Gamma\left( \frac {A+1}{a+1} \right) }{a+1}
  v^{-\frac {A+1}{a+1} } +
  \sum_{k=0}^M  \frac{(-1)^k}{k!} \zeta(-A-k(a+1)) v^{k} 
  + \frac 1{2\pi i} \int_{-M-\frac 12-i\infty}^{-M-\frac 12+i\infty}
  \Gamma(s) \zeta((a+1)s - A) v^{-s}\, ds \\
  &=& \frac {\Gamma\left( \frac {A+1}{a+1} \right) }{a+1} v^{-\frac
    {A+1}{a+1} } +
  \sum_{k=0}^M  \frac{(-1)^k}{k!} \zeta(-A-k(a+1)) v^{k} 
  + O\left( v^{M+\frac 12} \right),
\end{eqnarray*}
which proves the first part of the lemma.

If $(A+1)/(a+1) = -k_0$ for some integer $k_0\ge 0$ then $\Gamma(s)$
and $\zeta((a+1)s - A)$ create a double pole at $s = -k_0$ with
residue
\[
{\rm Res}\left( \Gamma(s) \zeta((a+1)s - A) v^{-s}, -k_0 \right) = 
\frac{(-1)^{k_0}}{(a+1)k_0!}\left( 1 + (H_{k_0} + a\gamma)
  \log \frac 1v \right) v^{k_0}
\]
of the resulting function.  This explains the difference from the first
case and completes the proof of the lemma.
\end{Proof}

We also need representations for finite sums of powers of integers
that can be deduced from the Euler--McLaurin summation formula, see e.g. 
\cite[Chapter 9]{GKP}.
\begin{lemma}\label{sums}
  We have the following representations or asymptotic series
  representation, resp., for the sums $\sum_{j=1}^n j^{a}$, $a>0$:
\begin{enumerate}[label=\emph{\alph*)}]
 \item If $a$ is a non-negative integer,
 \begin{equation}\label{aint}
   \sum_{j=1}^n j^{a} = \frac{n^{a+1}}{a+1} + \frac{n^a}2 + \sum_{k =
     1}^{\lfloor a/2 \rfloor } \frac{B_{2k}}{a+1} {a+1 \choose 2k}
   n^{a+1-2k},
\end{equation}
where $B_{2k}$ are the Bernoulli numbers.
\item If $a$ is a real number different from the non-negative
  integers, we have the asymptotic series expansion
\begin{equation}\label{anonint}
  \sum_{j=1}^n j^{a} \sim \zeta(-a) + \frac{n^{a+1}}{a+1} +
  \frac{n^a}2 + \sum_{k\ge 1} \frac{B_{2k}}{a+1} {a+1 \choose 2k}
  n^{a+1-2k}.
\end{equation}
\end{enumerate}
\end{lemma}

We are now ready to prove Theorem~\ref{Th1}.
Let us write $\K_{a,b}(x) = \sum_{i=1}^\infty P_i(x)$, where $P_i(x) :=
i^{b+1} \prod_{j=1}^i \frac{1}{1+ j^a/x}$, and note that, as $x\to \infty$,
\[
P_i(x) \sim i^{b+1} \exp\left( -\sum_{j=1}^i \log \left(1 + \frac{j^a}x
\right) \right).
\]
First of all we prove that $P_i(x)$ do not contribute 
significantly to $\K_{a,b}(x)$ if $i > x^{3/(3a+2)}$.

\begin{lemma}\label{bign}
We have  
\[
\sum_{i > x^{3/(3a+2)}}  P_i(x) \leq K e^{- c x^{1/(3a+1)}}
\]
for some constants $c,K> 0$.
\end{lemma}

\begin{Proof}
First, we assume that $i^a > \alpha x$, where $\alpha$ will be chosen
later. Then we have
\[
P_i(x) = i^{b+1} \prod_{j=1}^i \frac{1}{1+ j^a/x} \le  \prod_{j=1}^i
i^{1+b} \frac {x}{j^a} = i^{1+b} \frac{x^i}{(i!)^a} 
\]
Since $i! \ge (i/e)^i$ and $i^a > \alpha x$ we thus obtain
\[
  P_i(x) \le i^{b+1}\left( \frac {e^a x}{i^a} \right)^i \le
  \alpha^{-i}e^{ai+(1+b)i}.
\]
Consequently, if we choose $\alpha= e^{2a+(1+b)}$, we have that $ P_i(x)
\leq e^{-ai}$ and hence 
\[
\sum_{i:\, i^a > A x} P_i(x) \le \sum_{i:\, i^a > \alpha  x} e^{-ai}
\le K_1 e^{-c_1 x^{1/a}} \le K_1 e^{- c_2 x^{1/(3a+1)}}
\]
for some constants $c_2, K_1 > 0$.

We now assume that $x^{3/(3a+2)} < i \le \alpha^{1/a}x^{1/a}$, with
$\alpha$ chosen as above. In this case we have that $j^a/x \le \alpha$ so that
there exists a constant $c_3> 0$ such that
\[
\log\left(1 + \frac{j^a}x \right) \ge c_3 \frac {j^a} x 
\]
for all $j\le i$. Consequently there exists a constant $c_4> 0$ such that
\begin{eqnarray*}  
P_i(x) &~=~& i^{b+1}\exp\left( - \sum_{j=1}^i \log \left(1 +
    \frac{j^a}x \right)\right)\\
&~\le~& i^{b+1}\exp\left( - c_3 \frac 1x \sum_{j=1}^i j^a \right)
\le i^{b+1} \exp\left( - c_4 \frac{i^{a+1}}x \right).
\end{eqnarray*}
Note that for every real $\rho$, every $\sigma >0$ and
$\kappa_1, \kappa_2$ such that $0 <\kappa_1< \kappa_2$ there is a
constant $D$ depending on these four numbers such that for all $y>0$
\[
  y^\rho e^{-\kappa^2 y^\sigma} \leq D e^{-\kappa_1 y^\sigma}.
\]
Hence there are constants $K_2>0$ and $c_5, c_6>0$ such that 
\begin{eqnarray*}
  \sum_{x^{3/(3a+2)} < i \le \alpha^{1/a} x^{1/a}}   P_i(x) & \le &
  \alpha^{(b+1)/a} x^{(b+1)/a}
  \sum_{x^{3/(3a+2)} < i \le \alpha^{1/a} x^{1/a}} \exp\left(
    - c_4 \frac{i^{a+1}}x \right) \\
   &\le & \alpha^{(b+1)/a} x^{(b+1)/a} e^{- c_5 x^{1/(3a+1)}} \leq
  K_2 e^{- c_6 x^{1/(3a+1)}}.
\end{eqnarray*}
Now pick $c= \min\{ c_2, c_6\}$ and $K= \max\{ K_1, K_2\}$ to complete
the proof of the lemma.
\end{Proof}

Thus, it remains to consider $i$ with $i \le x^{3/(3a+2)}$. In this
case we certainly have $j^a/x \to 0$ as $x \rightarrow \infty$, so
we can use the Taylor expansion of $\log(1+z)$ to proceed
further.  From this Taylor series expansion, it follows that for every
$L\ge 1$ we have uniformly for $j\le i$
\[
  \log\left(1 + j^a/x \right) = \sum_{\ell=1}^{L-1} (-1)^{\ell-1}
  \frac 1{\ell} \frac{j^{\ell a}}{x^\ell} + O\left( \frac{j^{L
        a}}{x^L} \right),
\]
and consequently
\begin{equation}\label{eqrepsumlog}
\sum_{j=1}^i \log\left(1 + j^a/x \right) = \sum_{\ell=1}^{L-1}
(-1)^{\ell-1} \frac 1{\ell} \frac 1{x^\ell} \sum_{j=1}^i j^{\ell a} +
O\left( \frac{i^{aL+1}}{x^L} \right).
\end{equation}
In order to handle these terms we will use Lemma \ref{sums}.

With the help of the representation (\ref{eqrepsumlog}) and
Lemma~\ref{sums} we see that
$\sum_{j=1}^i \log\left(1 + j^{a\ell}/x \right)$ is dominated by
$i^{a\ell+1}/((a\ell+1)x)$ followed by smaller order terms. Here we
have to distinguish between $\ell = 1$ and $\ell \ge 2$. For $\ell =1$
the dominating term $i^{a+1}/((a+1)x)$ is unbounded if
$i \le x^{3/(3a+2)}$ whereas the next order term $i^a/x$ (and all
following terms) are bounded (in order) by $x^{-2/(3a+2)}$. So all of
them go to zero if $x\to\infty$.  It $\ell \ge 2$ then the dominant
terms $i^{a\ell+1}/((a\ell+1)x)$ (and, thus, all following terms) will
go to zero, too.  They are bounded (in order) by
$x^{(3-2\ell)/(3a+2)} \le x^{-1/(3a+1)}$.

Summing up, we obtain for $i \le x^{3/(3a+2)}$
\[
  \sum_{j=1}^i \log\left(1 + j^a/x \right) \sim \frac{i^{a+1}}{(a+1)
    x} + \tilde S,
\]
where $\tilde S$ collects terms of the form $\log \frac 1x$ that go to zero.
Hence by using the Taylor series of the exponential function we have
\[
P_i(x) \sim i^{b+1} \exp\left( -\frac{i^{a+1}}{(a+1) x} \right) \left(
1 + \tilde S + \frac 12 {\tilde S}^2 + \cdots \right)
\]
which leads again to an asymptotic series representation for $P_i(x)$
of the form
\begin{equation}\label{eqPiexp}
  P_i(x) \sim  i^{b+1} \exp\left( -\frac{i^{a+1}}{(a+1) x} \right) \left( 1 + \tilde T \right),
\end{equation}
where $\tilde T$ collects terms of the form $\hbox{const}\cdot
i^A/x^B$ (with real $A$ and integer $B\ge 1$) that go to zero if $i
\le x^{3/(3a+2)}$.

This discussion shows that we are finally led to consider sums of the
form
\[
\sum_{1\le i \le x^{3/(3a+2)}}  \exp\left(- \frac{i^{a+1}}{(a+1) x} \right) i^A.
\]
Since the sum of the missing terms can be estimated by
\[
\sum_{i > x^{3/(3a+2)}} \exp\left( \frac{i^{a+1}}{(a+1) x} \right) i^A
\le e^{-c_7 x^{1/(3a+1)}}
\]
for some constant
$c_7>0$, it is sufficient to consider infinite sums of the form
analysed in Lemma~\ref{LemmaR}.

We are now ready to complete the proof of Theorem~\ref{Th1}.

Since (\ref{eqPiexp}) are asymptotic series for $P_i(x)$, it follows
that we can consider them always as finite sums plus an error term of
the same form. Thus we can sum over them, at least for $i\le
x^{3/(3a+2)}$. However, by Lemma~\ref{bign} we can extend this
summation over all $i\ge 1$ since the resulting error is negligible.

Considering the terms in $\tilde{T}$, we see that the asymptotic
representation is different for the case of $a$ being a positive
integer, and for non-integer real positive $a$.

In the first case, by applying the (\ref{aint}) of Lemma~\ref{sums}
and observing that we just get positive powers of $n$ in the
representation of the sums $\sum_{j=1}^i j^{\ell a}$, $1\le \ell < L$,
the asymptotic series expansion (\ref{eqPiexp}) of $P_i(x)$, $i\le
x^{3/(3a+2)}$ can be written in the form
\[
  P_i(x)  \sim i^{b+1} \exp\left( -\frac{i^{a+1}}{(a+1) x} \right) \left( 1 +
    \sum_{B=1}^\infty \sum_{A=1}^{aB + \lfloor B/2 \rfloor} \tilde
    c_{A,B} \frac{ i^{A} }{x^{B}} \right),
\]
where $A,B$ are now integers and $\tilde c_{A,B}$ are real constants. 

This means that we also get an asymptotic series representation of
$\K_{a,b}(x)$ of the form
\begin{equation}\label{Kaint}
  \K_{a,b}(x) \sim R(a,1+b;1/((a+1)x)) + \sum_{B=1}^\infty 
\sum_{A=1}^{A=aB +\lfloor B/2 \rfloor} \tilde c_{A,B} \frac{R(a,1+b+A;1/((a+1)x))}{x^B}.
\end{equation}

For non-integer $a$ we we can proceed in the same way as in the
integer case. There are, however, some differences in the course of
the computations. First of all the sums $\sum_{j=1}^i j^{\ell a}$ do
not have an explicit representations. By (\ref{anonint}) of
Lemma~\ref{sums}, we obtain an asymptotic series expansion that
contains also negative powers of $i$, namely $i^{a+1-k}$ for any $k\ge
0$. This leads to an asymptotic series expansion for $P_i(x)$ of the
form
\[
  P_i(x) \sim i^{b+1}\exp\left( -\frac{i^{a+1}}{(a+1) x} \right) \left( 1 +
    \sum_{B=1}^\infty \, \sum_{A=-\infty}^{A=\lfloor B/2 \rfloor}
    c_{A,B} \frac{i^{A+aB}    }{x^{B}} \right),
\]
where $c_{A,B}$ are again real constants and the sum ranges over all
(even negative) integers $A\le \lfloor B/2 \rfloor$.  In completely
the same way as above we get from that an asymptotic series expansion
for $\K_{a,b}(x)$:
\begin{equation}\label{Kanonint}
  \K_{a,b}(x)  \sim R(a,1+b;1/((a+1)x)) + \sum_{B=1}^\infty\, \sum_{A=
    -\infty}^{A=\lfloor B/2\rfloor} c_{A,B} \frac{R(a,1+b+aB+A;1/((a+1)x))}{x^B}.
\end{equation}

Now we are going to use the information in Lemma~\ref{LemmaR} to
understand the leading terms and the order of the remainder in
(\ref{Kaint}) and (\ref{Kanonint}).  

In both cases, of integer and non-integer $a$, we have that for
any fixed integer $B$, the expression ${\displaystyle
  \frac{b+2+A}{a+1}-B}$ (for integer $a$) and ${\displaystyle
  \frac{b+2+aB +A}{a+1}-B}$ (for non-integer $a$) are maximised by
taking the largest allowable $A$ in these two cases to give
\[
\frac{b+2}{a+1}+\frac{1}{a+1}\left(\lfloor B/2\rfloor -B \right),
\]
and this is maximised by picking $B=1,2$ to give $(b+1)/(a+1)$.
Hence if $b>-1$, the asymptotic series (\ref{Kaint}) and
(\ref{Kanonint}) give, using (\ref{Ranonint})
\[
\K_{a,b}(x) \sim \frac{\Gamma\left( \frac{b+2}{a+1} \right)}{
(a+1)^{1-(b+2)/(a+1)}} x^{\frac{b+2}{a+1}} + O \left( x^{\frac{b+1}{a+1}}
\right). 
\]

If $b \leq -1$, we get 
\[
\K_{a,b}(x) \sim \frac{\Gamma\left( \frac{b+2}{a+1} \right)}{
(a+1)^{1-(b+2)/(a+1)}} x^{\frac{b+2}{a+1}} + O \left( 1 \right),  
\]
with the contributions to the $O(1)$ term coming from
$R(a,1+b;1/((a+1)x))$ term if $b<-1$. If $b=-1$, the contributions to
the $O(1)$ term come again from the $R(a,1+b;1/((a+1)x))$ term, and,
if $a$ is an integer, from the $B=1, A=1$ term in (\ref{Kaint}); if
$a$ is not an integer, the additional contributions come from the
$B=1, A=0$ and the $B=2, A=1$ terms in (\ref{Kanonint}). 

\medskip

{\bf Remarks:} 1. Much more can be said using Lemma~\ref{LemmaR}. For
example, if $b=-2$ the leading term of the asymptotic expansion of
$K_{a,b}(x)$ will be of order $\log x$. 

2. The case $b=-1$, $a=1$ can be solved explicitly. There we have
\[
\K_{1,-1}(x) \sim \sqrt{\frac{\pi}{2}} x^{1/2} - \frac23 +
\frac{\sqrt{2\pi}}{24} x^{-1/2} + O(1/x). 
\]

3. It is easy to show that 
\[
  K_{a,b}\left( \frac{x}{2} \right) \leq \H_{a,b}(x) \leq \K_{a,b}(x).
\]
Another case that can be solved explicitly is the case of $\H_{1,0}(x)$,
\[
  \H_{1,0}(x) = \sum_{i=1}^\infty i \prod_{j=1}^i \frac{x}{x+j+1}.
\]  
In this case MAPLE can compute the series and its asymptotics. We have
\[
  \H_{1,0}(x) = x -e^x(\Gamma(x+2)-\Gamma(x+2,x))x^{-x-1}, 
\]
where
\[
  \Gamma (x+2,x) = \int_x^\infty e^{-t} t^{x+1} \, dt
\]  
is the incomplete $\Gamma$-function. The asymptotic expansion as
$x \rightarrow \infty$ is
\[
 \H_{1,0}(x) \sim x - \sqrt{\frac{\pi}{2}} x^{1/2} + \frac53 + O(x^{-1/2}). 
\]

\section{Further Remarks and Conclusions} 

Corollary~\ref{Th2} is clearly non-optimal. Numerical evidence shows
that if $b \geq a-1$ the equilibrium is unique for all $\alpha/r$ and
there are at most two equilibria if $b < a-1$. To prove such a result
one might try to use the machinery of Pinelis \cite{Pinelis}. The result
for $b \geq a-1$ would follow if one could prove that $\H(x)$ is
convex and for $b < a-1$, by showing $\H(x)$ is concave.  However,
these results are difficult to obtain even for $\K_{a,b}(x)$. It does
seem that for each fixed $x$, the second derivative of $\K_{a,b}(x)$ is
an increasing function of $b$.

In the case of $\H_{1,0}(x)$ considered above, numerically $\H_{1,0}$
appears to be convex. Hence if we could prove that, and monotonicity
of the second derivative in $b$ for $\H_{a,b}(x)$ for every fixed $x$,
the desired results for $\H_{a,b}(x)$ would follow by the argument of
Pinelis \cite{Pinelis}.

If $b<a-1$, a possible strategy for proving, for example, that
$\K_{a,b}(x)$ is concave is to consider partial sums of the infinite
sum. For $b \geq a-1$ this strategy also works but is more interesting
because $x/(x+1)$ is concave. In that case, if $S_n(x)$ is the partial
sum, numerics indicate that $S_n''(x)$ is positive for $0<x<x_n$ and
$x_n \rightarrow \infty$ s $n \rightarrow \infty$.

Note that for all integer values of $a, b$ MAPLE can compute
$\F_{a,b}(x)$ in terms of hypergeometric functions. This, however,
does not seem very useful.

From the biological point of view, $b$ measures the efficiency of the
muco-ciliatory escalator, while $a$ measures its inefficiency due to
release of quartz in the lungs by macrophages with supercritical
load. Our results show that the ratio $(b+2)/(a+1)$ is crucial in
establishing whether the system can deal with the quartz load; if it
is less or equal to 1, there is a deposition rate $\alpha$ that will
overwhelm it, no matter what $r$ is.

In summary, we have completed the model of \cite{Tran} by including an
equation for the evolution of quartz concentration. The resulting
mathematical object is a challenging system of coagulation--death
equations that requires non-trivial asymptotic ideas in the discussion
of the structure of equilibria. Of course the analysis in the paper is
only part of the necessary mathematical work; one also needs to
establish global existence (using finite-dimensional truncations or
methods of semigroup theory) and stabilisation to equilibria (for
example, by exhibiting a suitable Lyapunov function). 

\subsection*{Acknowledgements}

FPdC and MG would like to thank P. Freitas and D. Pritchard for
valuable discussions. FPdC acknowledges financial support provided by
the University of Strathclyde David Anderson Research Professorship,
and all the authors are grateful to D. Zeilberger for helping to set
up the collaboration.

\end{document}